\documentclass[11p,reqno]{amsart}

\usepackage[T1]{fontenc}
\usepackage{graphicx}
\usepackage{amssymb,amsthm,amsmath,mathrsfs,bm,braket,marginnote}
\usepackage{enumerate}
\usepackage{appendix}
\usepackage[colorlinks=true, pdfstartview=FitV, linkcolor=blue, citecolor=blue, urlcolor=blue]{hyperref}
\usepackage{multirow}

\usepackage{pgf}
\usepackage{pgfplots}
\usepackage{tikz}
\usetikzlibrary{arrows,calc}
\usepackage{verbatim}
\usetikzlibrary{decorations.pathreplacing,decorations.pathmorphing}
\usepackage[numbers,sort&compress]{natbib}
\usepackage{dsfont}

\numberwithin{equation}{section}
\newtheorem{theorem}{Theorem}[section]

\newtheorem{remark}[theorem]{Remark}

\newtheorem{proposition}[theorem]{Proposition}

 \reversemarginpar
 
 \newcommand{\Pf}{\text{Pf}}
\newcommand{\dq}{d_qx}
\newcommand{\al}{{\alpha}}
\newcommand{\be}{{\beta}}

\begin{document}

\title[Evaluation of certain Pfaffians ]{Evaluations of certain Catalan-Hankel Pfaffians via classical skew orthogonal polynomials}

\author{Bo-Jian Shen}
\address{School of Mathematical Sciences, Shanghai Jiaotong University, People's Republic of China.}
\email{JOHN-EINSTEIN@sjtu.edu.cn}
\author{Shi-Hao Li}
\address{Department of Mathematics, Sichuan University, Chengdu, 610064, China}
\email{lishihao@lsec.cc.ac.cn}
\author{Guo-Fu Yu}
\address{School of Mathematical Sciences, Shanghai Jiaotong University, People's Republic of China.}
\email{gfyu@sjtu.edu.cn}

\date{}

\dedicatory{}


\begin{abstract}
{This paper is to evaluate certain Catalan-Hankel Pfaffians by the theory of skew orthogonal polynomials. Due to different kinds of hypergeometric orthogonal polynomials underlying the Askey scheme, we explicitly construct the classical skew orthogonal polynomials and then give different examples of Catalan-Hankel Pfaffians with continuous and $q$-moment sequences.}
\end{abstract}

\maketitle

\section{Introduction}
Hankel determinants as a specific determinant has attracted much attention from researchers in many different subjects (see, for example \cite{bacher02,han16,vein} and references therein). 
It is well known that if $\{\mu_m\}_{m\geq0}$ is the moment sequence taking the form
\begin{align*}
\mu_n=\int_{\mathbb{R}}x^n \omega(x)dx,
\end{align*}
then the Hankel determinant $\det(\mu_{i+j})_{i,j=0}^{n-1}$ has a nice integral formula due to the Andrei\'ef formula \cite{forrester19}
\begin{align*}
\det(\mu_{i+j})_{i,j=0}^{n-1}=\frac{1}{n!}\int_{\mathbb{R}^n}\prod_{1\leq i<j\leq n}|x_i-x_j|^2\prod_{i=1}^n \omega(x_i)dx_i.
\end{align*}
Moreover, if $\omega(x)$ is a classical weight, e.g. the Gaussian weight, Laguerre weight, Jacobi weight and circular Jacobi weight, then Hankel determinant with classical moments can be explicitly evaluated due to the Selberg integral \cite[Chapter 4]{forrester10}. Moreover, $q$-versions of Selberg integrals which relate to the moment sequence with discrete measure play important roles in modern mathematical physics as well \cite{forrester08}.

In recent years, the Pfaffian version of Hankel determinant, called the Catalan-Hankel Pfaffian, was proposed due to its potential applications in the theory of combinatorics \cite{ishikawa13}. The Pfaffian version of Andrei\'ef formula, which is called the de Bruijn formula, has been applied into the evaluation of certain Catalan-Hankel Pfaffian with the help of the Selberg integral and its $q$-version \cite{ishikawa13,ishikawa20}. Catalan-Hankel Pfaffian has the form
\begin{align}\label{cpf}
\Pf\left[(j-i)\mu_{i+j-1}\right]_{i,j=0}^{2N-1}
\end{align}
and its $q$-version
\begin{align}\label{dpf}
\Pf\left[
([j]_q-[i]_q)\mu_{i+j-1}
\right]_{i,j=0}^{2N-1},\quad [j]_q=\frac{1-q^j}{1-q},
\end{align}
where $\{\mu_n\}_{n\geq0}$ is a moment sequence related to discrete $q$-measure. The specific little $q$-Jacobi case were considered in \cite{ishikawa13} and some recent developments like the Al-Salam $\&$ Carlitz I case was given in \cite{ishikawa20}. 

Regarding evaluations of ($q$-)Catalan-Hankel Pfaffian, one application is to give formula to the weighted enumeration of some specific plane partitions. For example, in \cite{ishikawa13}, the authors considered the little $q$-Jacobi case and enumerated a special family of shifted reverse plane partitions with weights that resemble the weight in the inner product of little $q$-Jacobi polynomials. 
The importance of these evaluations lies in the random matrix theory as well.
In an earlier work of Mehta and Wang \cite{mehta00}, they gave an evaluation of the Hankel Pfaffian
\begin{align}\label{mehta}
\Pf[(j-i)\Gamma(i+j+b)]_{i,j=0}^{2N-1},
\end{align}
which was related to the skew orthogonal Laguerre polynomials given in \cite{adler00}.
Moreover, the shifted Catalan-Hankel Pfaffian 
\begin{align*}
\Pf\left[
(j-i)\mu_{i+j+r}\right]_{i,j=0}^{2N-1}
\quad \text{or}\quad
\Pf\left[
([j]_q-[i]_q)\mu_{i+j+r}
\right]_{i,j=0}^{2N-1}
\end{align*}
is closely related to the adjacent families of skew-orthogonal polynomials considered in \cite{li20}. These evaluations may give more hints in exploring novel examples in integrable systems and random matrices.

As mentioned, one way to evaluate these Pfaffians is based on the de Bruijn formula and Selberg integral \cite{ishikawa13,ishikawa20}. In this paper, we will investigate another way, namely by using the theory of classical ($q$-)skew orthogonal polynomials, to make these evaluations and demonstrate the effectiveness---by simply considering different classsical weights underlying Askey scheme, we can get different examples of ($q$-)Catalan-Hankel Pfaffians.
Since the Hankel determinant is closely related to the theory of classical orthogonal polynomials, then 
if we know normalisation constants of these classical orthogonal polynomials, the evaluations of Hankel determinant can be made. It is natural to ask whether we can apply the theory of classical ($q$-)skew orthogonal polynomials into evaluations of the certain ($q$-)Catalan-Hankel Pfaffians.
The answer is affirmative. If we know the normalisation factors of skew orthogonal polynomials, then these Pfaffians can be explicitly given by those normalisation factors. Due to different expressions in discrete and continuous cases  (c.f. \eqref{cpf} and \eqref{dpf}), we discuss them separately.
Continuous cases including Hermite-type, Laguerre-type and Jacobi-type skew orthogonal polynomials were firstly constructed by Adler et al \cite{adler00} and we give a brief review in section \ref{sec2}. Some evaluations of Catalan-Hankel Pfaffians related to these weights including formula \eqref{mehta} will be demonstrated. Moreover, we consider a Catalan-Hankel Pfaffian with Cauchy weight  in that part. Evaluations of $q$-Catalan-Hankel Pfaffians will be based on the theory of classical $q$-skew orthogonal polynomials. With the help of discrete Pearson relation, we construct different examples of classical ($q$-)skew orthogonal polynomials, thus obtaining different kinds of skew orthogonal polynomials and normalisation factors with respect to different classical $q$-weights, which enlarge the results given by \cite{forrester20,ishikawa20}.

\section{Continuous measure}\label{sec2}
Let's consider the skew inner product related to \eqref{cpf}
\begin{align}\label{sip}
\langle \phi(x),\psi(x)\rangle_{4,\omega}=\frac{1}{2}\int_{\mathbb{R}}[\phi(x)\psi'(x)-\phi'(x)\psi(x)]\omega(x)dx,
\end{align}
then the skew moments are given by 
\begin{align*}
m_{i,j}:=\langle x^i,x^j\rangle_{4,\omega}=\frac{1}{2}(j-i)\int_{\mathbb{R}}x^{i+j-1}\omega(x)dx.
\end{align*}
The skew LU decomposition (or so-called skew Borel decomposition \cite{adler99}) of the moment matrix $(m_{i,j})_{i,j\in\mathbb{N}}$ could give rise to 
the monic skew orthogonal polynomials $\{Q_j(x)\}_{j\in\mathbb{N}}$ satisfying the skew orthogonal relation
\begin{align}\label{scd}
\langle Q_{2n}(x),Q_{2m+1}(x)\rangle_{4,\omega}=u_n\delta_{n,m},\quad \langle Q_{2n}(x),Q_{2m}(x)\rangle_{4,\omega}=\langle Q_{2n+1}(x),Q_{2m+1}(x)\rangle_{4,\omega}=0
\end{align}
for certain $u_n>0$. 
Moreover, these skew orthogonal polynomials have the following Pfaffian expressions \cite{chang18}
\begin{align*}
Q_{2n}(x)=\frac{1}{\tau_{2n}}\Pf(0,\cdots,2n,x),\quad Q_{2n+1}(x)=\frac{1}{\tau_{2n}}\Pf(0,\cdots,2n-1,2n+1,x)
\end{align*}
with $\tau_{2n}=\Pf(0,\cdots,2n-1)$ and the Pfaffian elements are given by $\Pf(i,j)=m_{i,j}$ and $\Pf(i,x)=x^i$. By putting $m_{i,j}$ into the expression of $\tau_{2n}$, we have
\begin{align}\label{hs}
\tau_{2n}=\frac{1}{2^n}\Pf\left[(j-i)\mu_{i+j-1}\right]_{i,j=0}^{2n-1}, \quad \mu_n=\int_{\mathbb{R}}x^n\omega(x)dx.
\end{align}
Interestingly, $u_n$ in the skew orthogonal condition \eqref{scd} is the ratio $\tau_{2n+2}/\tau_{2n}$, and 
therefore, if we can explicitly compute $\{u_n\}_{n\in\mathbb{N}}$ in the skew orthogonal condition, then we directly have
$
\tau_{2N}=\prod_{i=0}^{N-1}u_i.
$
In fact, normalisation factors $\{u_n\}_{n\in\mathbb{N}}$ can be explicitly computed when $\{Q_n(x)\}_{n\in\mathbb{N}}$ are {\bf{classical}} skew orthogonal polynomials. For details, please refer to \cite{adler00} and we give a brief review here.

Starting from the inner product
\begin{align*}
\langle \phi(x),\psi(x)\rangle_{2,\rho}=\int_{\mathbb{R}}\phi(x)\psi(x)\rho(x)dx,
\end{align*}
one can construct a family of \textbf{monic} orthogonal polynomials $\{p_j(x)\}_{j\in\mathbb{N}}$ satisfying the orthogonal relation
\begin{align*}
\langle p_j(x),p_k(x)\rangle_{2,\rho}=h_j\delta_{j,k}.
\end{align*}
We call these orthogonal polynomials classical if the weight function $\rho(x)$ satisfies the following Pearson equation
\begin{align*}
\frac{\rho'(x)}{\rho(x)}=-\frac{g(x)}{f(x)} \quad\text{with deg $f(x)\leq 2$ and deg $g(x)\leq 1$}.
\end{align*}
From the relation, one can construct an operator 
\begin{align}\label{co}
\mathcal{A}=f\partial_x+\frac{f'-g}{2},
\end{align}
such that
\begin{align*}
\langle \phi(x),\mathcal{A}\psi(x)\rangle_{2,\rho}=\langle \phi(x),\psi(x)\rangle_{4,\omega},\quad \omega(x)=\rho(x)f(x),
\end{align*}
where the skew inner product $\langle \cdot,\cdot\rangle_{4,\omega}$ was given in \eqref{sip}. This is the key formula to establish the relation between classical orthogonal polynomials and classical skew orthogonal polynomials. Moreover, the normalisation factor $u_n$ could be computed via
\begin{align}\label{ana}
\mathcal{A}p_k(x)=-\frac{c_k}{h_{k+1}}p_{k+1}(x)+\frac{c_{k-1}}{h_{k-1}}p_{k-1}(x),\quad u_{n}=c_{2n},
\end{align}
with $h_k$ the normalisation factor of orthogonal polynomials. Therefore, to obtain normalisation factors of skew orthogonal polynomials, our attention should be paid to the computation of $c_k$ in the above equation. The fastest method to compute $c_k$ is to compare the coefficients of $x^{k+1}$ on both sides and the following are the examples of the continuous classical weights including Hermite, Laguerre, Jacobi and Cauchy weights.

From weight functions 
\begin{align*}
\rho(x)=\left\{
\begin{array}{ll}
e^{-x^2},& \text{Hermite},\\
x^ae^{-x},& \text{Laguerre},\\
x^a(1-x)^b,& \text{Jacobi},\\
(1+x^2)^{-a},& \text{Cauchy},
\end{array}
\right.
\end{align*}
one can obtain the Pearson pair
\begin{align*}
(f,g)=\left\{
\begin{array}{ll}
(1,2x),&\text{Hermite},\\
(x,x-a),&\text{Laguerre},\\
(x(1-x),(a+b)x-a),&\text{Jacobi},\\
(1+x^2,2ax),&\text{Cauchy}.
\end{array}
\right.
\end{align*}
This table was given by \cite[eq. (5.58)]{forrester10} and it should be remarked that we use the weight function of Jacobi as $x^a(1-x)^b$ supported in $[0,1]$ to make the moments easily written down. 
There are two quantities obtained from this table. One is the Hankel sequences given in \eqref{hs}.
From $\omega(x)=f(x)\rho(x)$, one knows that the weights of classical skew orthogonal polynomials are
\begin{align*}
\omega^{(H)}(x)=e^{-x^2},\quad \omega^{(L)}(x)=x^{a+1}e^{-x},\quad \omega^{(J)}(x)=x^{a+1}(1-x)^{b+1},\quad \omega^{(C)}(x)=(1+x^2)^{-a+1}.
\end{align*}
Therefore, we have the following moments form
\begin{align*}
&\mu_{n}^{(H)}=\frac{1+(-1)^{n}}{2}\Gamma\left(
\frac{n+1}{2}
\right),\quad \mu_{n}^{(L)}=\Gamma(n+a+2),\\
&\mu_{n}^{(J)}=\frac{\Gamma(n+a+2)\Gamma(b+2)}{\Gamma(n+a+b+4)},\quad\,\,\, \mu_{n}^{(C)}=\frac{1+(-1)^{n}}{2}\frac{\Gamma((n+1)/2)\Gamma(a-2-(n+1)/2)}{\Gamma(a-1)}.
\end{align*}
On the other hand, we can obtain 
\begin{align*}
c_k=\left\{\begin{array}{ll}
h_{k+1}^{(H)},&\text{Hermite},\\
h_{k+1}^{(L)}/2,&\text{Laguerre},\\
(k+1+(a+b)/2)h_{k+1}^{(J)},&\text{Jacobi},\\
(a-1-k)h_{k+1}^{(C)},&\text{Cauchy}
\end{array}
\right.
\end{align*}
where $\{h_{k}\}_{k\in\mathbb{N}}$ are normalisation factors with respect to different weights \cite[Chap. 5]{forrester10}
\begin{align*}
&h_k^{(H)}=\pi^{1/2}2^{-k}k!,\quad h_k^{(L)}=\Gamma(k+1)\Gamma(a+k+1),\\
&h_k^{(J)}=2^{a+b+1+2k}\frac{\Gamma(k+1)\Gamma(a+b+1+k)\Gamma(a+1+k)\Gamma(b+1+k)}{\Gamma(a+b+2k+1)\Gamma(a+b+2k+2)},\\
&h_k^{(C)}=\pi 2^{2k+2a+2}\frac{\Gamma(k+1)\Gamma(-2k-2a)\Gamma(-2a-2k-1)}{\Gamma(-2a-k)(\Gamma(-a-k))^2}.
\end{align*}

By using above results, we could state the following proposition.
\begin{proposition}
We have the following evaluations of certain Catalan-Hankel Pfaffians
\begin{align*}
&\Pf\left[(j-i)\mu_{i+j-1}^{(H)}\right]_{i,j=0}^{2N-1}=2^{-N(N-1)}(\sqrt{\pi})^N\prod_{i=0}^{N-1}\Gamma(2i+2),\\
&\Pf\left[(j-i)\mu_{i+j-1}^{(L)}\right]_{i,j=0}^{2N-1}=\prod_{i=0}^{N-1}\Gamma(2i+2)\Gamma(2i+a+2),\\
&\Pf\left[(j-i)\mu_{i+j-1}^{(J)}\right]_{i,j=0}^{2N-1}=4^N\frac{\Gamma(N+(a+b+2)/4)}{\Gamma((a+b+2)/4)}\prod_{i=0}^{N-1}\frac{\Gamma(a+2i+2)\Gamma(b+2i+2)\Gamma(2i+3)}{\Gamma(a+b+4i+3)},\\
&\Pf\left[(j-i)\mu_{i+j-1}^{(C)}\right]_{i,j=0}^{2N-1}=(-1)^N2^{2N^2-N(a-3)}\frac{\Gamma(N+(1-a)/2)}{\Gamma((1-a)/2)}\prod_{i=0}^{N-1}\frac{\Gamma(2a-4i-2)\Gamma(2i+3)}{(\Gamma(a-2i-1))^2}.
\end{align*}
\end{proposition}
\begin{remark}
The second formula is exactly equation \eqref{mehta} with the shift $a\to a-1$.
\end{remark}

\section{Discrete measure: $q$-case}
This part is devoted to the evaluations of $q$-Catalan Hankel Pfaffians given by formula \eqref{dpf}.
The $q$-case corresponds to a special discrete measure distributed on exponential lattices $ x(i)=q^i$ with $0<q<1$ and $i\in\mathbb{Z}$. 
By using the definition of Jackson's $q$-integral \footnote{The definition of $q$-integral on the interval $[0,\infty)$ is different from the one defined on $[0,1]$, see \cite{ismail05}. However, these two cases can be treated similarly and we just consider the former one here.}
\begin{align*}
	\int_{0}^{\infty}f(x)d_qx=(1-q)\sum_{s=-\infty}^{\infty}f(q^s)q^s,
\end{align*}
one can define the following inner product
\begin{align}\label{qip}
	\langle \phi(x),\psi(x)\rangle_{2,\rho}&:= \int_{0}^{\infty}\phi(x)\psi(x)\rho(x)d_qx=(1-q)\sum_{s\in \mathbb{Z}}\phi(q^s)\psi(q^s)\rho(q^s)q^s.
\end{align}
With this inner product, $q$-orthogonal polynomials  $\{p_n(x;q)\}_{n\in\mathbb{N}}$ are defined by the orthogonal relation
\begin{align}\label{qor}
\langle p_i(x;q),p_j(x;q)\rangle_{2,\rho}=h_n(q)\delta_{n,m}
\end{align}
with respect to the weight $\rho(x;q)$.
Moreover, if $\rho(x;q)$ is a classical weight, we call the corresponding orthogonal polynomials the classical $q$-orthogonal polynomials. Classical $q$-orthogonal polynomials include many interesting examples underlying the $q$-Askey scheme. For details, please refer to \cite{andrews00,ismail05,koekoek10}.
One important property of classical $q$-orthogonal polynomials is that the weight function satisfies an analogy of Pearson relation given by Nikiforov and Suslov \cite{nikiforov86}
\begin{align}\label{pearson}
\frac{\rho(qx)}{\rho(x)}=\frac{f(x)-q^{-\frac{1}{2}}(1-q)xg(x)}{f(qx)}	\quad\text{with deg $f(x)\leq 2$ and deg $g(x)\leq 1$}.
\end{align}
In the following, we will demonstrate how to connect $q$-inner product with $q$-skew inner product.

Let's define a $q$-analogy of the skew inner product \eqref{sip}
\begin{align}
	\langle \phi(x),\psi(x)\rangle_{4,\omega}=\int_{0}^{\infty}[\phi(x)D_q\psi(x)-D_q\phi(x)\psi(x)]\omega(x)d_qx\label{qsip}
\end{align}
with $q$-difference operator $$D_qf(x)=\frac{f(x)-f(qx)}{(1-q)x}.$$
Then by defining an operator $\mathcal{A}_q$ \cite{forrester20}
\begin{align}\label{qo}
\mathcal{A}_q=q^{-\frac{1}{2}}g(x)T_q+q^{-1}f(x)D_{q^{-1}}+f(x)D_q, \quad T_q f(x)=f(qx),
\end{align}
one can find such a connection formula
\begin{align}
\langle \phi(x;q),\mathcal{A}_q\psi(x;q)\rangle_{2,\rho}=\langle \phi(x;q),\psi(x;q)\rangle_{4,\omega},\quad \omega(x)=\rho(qx)f(qx).\label{oso}
\end{align}
On the other hand, from the skew inner product \eqref{qip}, one has the following ($q$-)skew moments
\begin{align}\label{sm}
m_{i,j}:=([j]_q-[i]_q)\int_{0}^{\infty}x^{i+j-1}\omega(x)d_qx, \quad [j]_q=\frac{1-q^j}{1-q}.
\end{align}
Similarly, it is known that a family of {\textbf{monic}} skew orthogonal polynomials $ \{Q_i(x;q)\}_{i\in\mathbb{N}} $ could be constructed from those moments if even-order moment matrices $\Pf(m_{i,j})_{i,j=0}^{2n-1}\neq0$ for all $n\in\mathbb{N}_+$ \cite{forrester20}. Furthermore, polynomials $\{Q_i(x;q)\}_{i\in\mathbb{N}}$ admit the following Pfaffian expressions 
\begin{align*}
&Q_{2n}(x;q)=\frac{1}{\tau_{2n}(q)}\Pf(0,\cdots,2n,x),\quad Q_{2n+1}(x;q)=\frac{1}{\tau_{2n}(q)}\Pf(0,\cdots,2n-1,2n+1,x)\\
&\tau_{2n}(q)=\Pf(0,\cdots,2n-1),\quad \Pf(i,j)=m_{i,j},\quad \Pf(i,x)=x^i,
\end{align*}
and they satisfy the following skew orthogonal relations
\begin{align*}
& \langle Q_{2n}(x;q),Q_{2m+1}(x;q)\rangle_{4,\omega}=\frac{\tau_{2n+2}(q)}{\tau_{2n}(q)}:=u_n(q)\delta_{n,m},\\
& \langle Q_{2n}(x;q),Q_{2m}(x;q)\rangle_{4,\omega}=\langle Q_{2n+1}(x;q),Q_{2m+1}(x;q)\rangle_{4,\omega}=0.
\end{align*}
Therefore, one knows that
\begin{align*}
\Pf\left[([j]_q-[i]_q)\int_{0}^{\infty}x^{i+j-1}\omega(x)d_qx\right]_{i,j=0}^{2N-1}=\Pf(0,\dots ,2N-1)=\prod_{i=0}^{N-1}u_i(q).
\end{align*}
Interestingly, there is a method to evaluate the value of $u_n(q)$ quickly by taking advantage of $q$-Pearson relation.

As an analogy of equation \eqref{ana}, there holds the formula
\begin{align}\label{eq:cn}
\mathcal{A}_q p_k(x;q)=-\frac{c_k(q)}{h_{k+1}(q)}p_{k+1}(x;q)+\frac{c_{k-1}(q)}{h_{k-1}(q)}p_{k-1}(x;q)
\end{align}
in the $q$-case due to the property of $\mathcal{A}_q$ \cite[eq. (4.28)]{forrester20}. Moreover, the quantity $h_k(q)$ is the normalisation factor of the orthogonal relation \eqref{qor} and $c_k(q)$ is closely related to $u_k(q)$ via the relation $u_k(q)=c_{2k}(q)$. Therefore, it is the point to compute $c_k(q)$ from above equation and then the exact value of $q$-Catalan Hankel Pfaffian could be obtained.
In the following part, we discuss several different cases including the examples in \cite{ishikawa13,ishikawa20}.

\subsection{Al-Salam $\&$ Carlitz I case}

The first example considered here is the Al-Salam $\&$ Carlitz I case with the weight function
\begin{align*}
\rho(x;q)=(qx,a^{-1}qx;q)_\infty,\quad a<0.
\end{align*}
It is well known that the Al-Salam $\&$ Carlitz polynomials $\{U_n^{(a)}(x;q)\}_{n\geq0}$ have the orthogonality
\begin{align*}
\int_a^1 U_m^{(a)}(x;q)U_n^{(a)}(x;q)\rho(x;q)d_qx=(-a)^n (1-q) (q;q)_n(q,a,a^{-1}q;q)_\infty q^{n\choose 2}\delta_{n,m}:={h}_n\delta_{n,m},
\end{align*}
and canonical moments were given by \cite{njionou}
\begin{align*}
	\mu_{n}=\int_{a}^{1}x^{n}\rho(x;q)d_qx=(1-q)\left(q, a, a^{-1} q ; q\right)_{\infty} \sum_{i=0}^{n}\left[\begin{array}{c}
	n \\
	i
	\end{array}\right]_{q} a^{i}
\end{align*}
Specifically, Al-Salam $\&$ Carlitz I polynomials have $q$-hypergeometric function expressions
\begin{align*}
U_n^{(a)}(x;q)=(-a)^n q^{n\choose 2} {_2}\phi_1\left(
\left.
\begin{array}{c}
q^{-n},x^{-1}\\
0
\end{array}\right|
q;\frac{qx}{a}
\right).
\end{align*}
As was shown in \cite[Sec. 4.4.1(1)]{forrester20}, the Pearson pair in this case is
\begin{align*}
(f,g)=\left(
x^2-(1+a)x+a,\,\frac{q^{1/2}}{1-q}(x-(1+a))
\right),
\end{align*}
and the coefficient $c_n$ in \eqref{eq:cn} can be explicitly computed as $c_n=-q^{-n}{h}_{n+1}/(1-q)$.
Thus the moment Pfaffian can be written as
\begin{align*}
	\Pf(m_{i,j})_{i,j=0}^{2N-1}&=\prod_{k=0}^{N-1}a^{2k+1}(q ; q)_{2k+1} q^{\binom{2k}{2}}\left(q, a, a^{-1} q ; q\right)_{\infty}\\&=a^{N^2}q^{\frac{1}{6}N(N-1)(4N-5)}\left(q, a, a^{-1} q ; q\right)_{\infty}^{N}\prod_{k=0}^{N-1}(q ; q)_{2k+1},
\end{align*}
where $ m_{i,j}=([j]_q-[i]_q)\int_{a}^{1}x^{i+j-1}\omega(x;q)\dq $. Since $ \omega(x;q)=f(qx;q)\rho(qx;q)=a\rho(x;q) $, we have
\begin{align*}
m_{i,j}=a([j]_q-[i]_q)\mu_{i+j-1}=(q^i-q^j)\left(q, a, a^{-1} q ; q\right)_{\infty} \sum_{k=0}^{i+j-1}\left[\begin{array}{c}
i+j-1 \\
k
\end{array}\right]_{q} a^{k+1}.
\end{align*}
Dividing both sides by $ \left[a(q, a, a^{-1} q ; q)_\infty\right]^{N} $ leads to
\begin{align*}
	\Pf\left((q^i-q^j)\sum_{k=0}^{i+j-1}\left[\begin{array}{c}
	i+j-1 \\
	k
	\end{array}\right]_{q} a^{k}\right)_{i,j=0}^{2N-1}=a^{N(N-1)}q^{\frac{1}{6}N(N-1)(4N-5)}\prod_{k=0}^{N-1}(q ; q)_{2k+1}.
\end{align*}

\begin{remark}
With $a=-1$, Al-Salam $\&$ Carlitz I polynomials reduce to the $q$-Hermite I polynomials \cite[Sec. 14.28]{koekoek10}, therefore the above mentioned method can be applied to the $q$-Hermite I case as well. The moments of q-Hermite I polynomials take the form \cite{njionou}
\begin{align*}
	\mu_{n}=(1-q)(q,-1,-q ; q)_{\infty} \frac{1+(-1)^{n}}{2}\left(q ; q^{2}\right)_{n / 2}.
\end{align*}
Thus we have the following evaluation of $q$-Catalan-Hankel Pfaffian 
\begin{align*}
	\Pf\left((q^i-q^j)\frac{1+(-1)^{i+j-1}}{2}(q;q^2)_{\frac{i+j-1}{2}}\right)_{i,j=0}^{2N-1}=(-1)^Nq^{\frac{1}{6}N(N-1)(4N-5)}\prod_{k=0}^{N-1}(q;q)_{2k+1}.
\end{align*}
\end{remark}

\subsection{Stieltjes-Wigert case}
The Stieltjes-Wigert polynomials are well studied in random matrix theory of the so-called Stieltjes-Wigert ensemble, which firstly appeared in the study of non-intersecting Brownian walkers, and subsequently in quantum many body systems etc.. For a detailed review, please refer to \cite{forrester20sw} and references therein.
The weight of Stieltjes-Wigert polynomials was first given by Stieltjes as an example of indeterminate moment problems \cite{stieljes} and further studied by Wigert \cite{wigert}. 
Usually there are several different expressions for the Stieltjes-Wigert's weight function \cite{chris03}. The original one corresponding to a continuous measure is
\begin{align*}
	w(x)=\frac{1}{\sqrt{\pi}} k x^{-k^{2} \log x}, \quad x>0
\end{align*}
with the moment $ \mu_{n}=\int_{0}^{\infty} x^{n} w(x) d x=e^{(n+1)^{2} / 4 k^{2}} $. If we set $ q=e^{-1/2k^2} $, then $\mu_n$ can be written as $ q^{-(n+1)^2/2} $. With this measure, Wigert found the following expression for Stieltjes-Wigert polynomials \cite{wigert}
\begin{align*}
	P_{n}(x)=(-1)^{n} \frac{q^{n / 2+1 / 4}}{\sqrt{(q ; q)_{n}}} \sum_{k=0}^{n}\left[\begin{array}{l}
	n \\
	k
	\end{array}\right]_{q}(-1)^{k} q^{k^{2}+k / 2} x^{k},
\end{align*}
with orthogonality
\begin{align*}
	\int_{0}^{\infty}P_{n}(x)P_{m}(x)w(x)dx=\delta_{m n}.
\end{align*}
In \cite{chihara70}, Chihara proposed a discrete weight on the exponential lattices admitting the form
\begin{align*}
\xi(x)=\left\{
\begin{array}{ll}
\frac{1}{\sqrt{q}M}  q^{n+n^2/2} & , x=q^n\\
0 & ,  x\neq q^n
\end{array} \right. ,\quad M=(-q\sqrt{q},-q^{-1/2},q;q)_\infty,\quad n\in \mathbb{Z}.
\end{align*}
Then the corresponding inner product becomes
\begin{align*}
	\langle p_n(x),p_m(x)\rangle_{2,\xi}=\frac{1}{\sqrt{q}M}\sum_{k=-\infty}^{\infty}p_n(q^k)p_m(q^k)q^{k+k^2/2}
\end{align*}

We can prove that this discrete weight is equivalent to the continuous one by showing that they have same moments. For completeness, we give a short proof to this fact.

One can easily check that
\begin{align*}
	\xi(qx)=q^{\frac{3}{2}}x\xi(x),\quad  x\in \mathbb{R}
\end{align*}
and therefore obtain a recurrence relation for the moments $ \{\mu_{n}\}_{n\in\mathbb{N}} $
\begin{align*}
	\mu_{n}&:=\langle x^n,1\rangle_{2,\xi}=\sum_{k=-\infty}^{\infty}q^{nk}\xi(q^k)=\sum_{k=-\infty}^{\infty}q^{n(k+1)}\xi(q^{k+1})\\
	&\quad=q^{n+\frac{3}{2}}\sum_{k=-\infty}^{\infty}q^{(n+1)k}\xi(q^k)=q^{n+\frac{3}{2}}\mu_{n+1}.
\end{align*}
As a result, we have $ \mu_n=q^{-n^2/2-n}\mu_0 $. By making use of the Jacobi triple product identity \cite[Thm 352]{hardy}
\begin{align*}
	\sum^{\infty}_{n=-\infty}(-1)^nq^{\binom{n}{2}}x^n=(x,q/x,q;q)_\infty,
\end{align*}
$ \mu_0 $ can be directly computed as
\begin{align*}
	\mu_0=\frac{1}{\sqrt{q}M}\sum^{\infty}_{k=-\infty}q^{k+k^2/2}=\frac{1}{\sqrt{q}}.
\end{align*}
Thus $ \mu_n=q^{-\frac{(n+1)^2}{2}} $, which is the same as that of the continuous measure.

For the present need, we would consider the discrete measure. Denote 
\begin{align*}
 \rho(x;q)=\frac{1}{\sqrt{q}M}x^{\frac{\ln x}{2\ln q}} 
 \end{align*}
  to be the weight function and define 
  \begin{align*}
  p_n(x)=\frac{\sqrt{(q;q)_n}}{q^{n^2+n+\frac{1}{4}}}P_n(x) 
  \end{align*}
  to be the monic Stieltjes-Wigert polynomials. We can rewrite the orthogonality in terms of $ \rho(x) $ and $ \{p_n(x)\}_{n\geq 0} $ as
\begin{align*}
	\int_{0}^{\infty}p_n(x)p_m(x)\rho(x)\dq=(1-q)\frac{(q;q)_n}{q^{2n^2+2n+\frac{1}{2}}}\delta_{m n}:=h_n\delta_{n,m}.
\end{align*}
By assuming $ f(x)=x $ and solving the Pearson equation \eqref{pearson},
we get $ g(x)={(q^2x-q^{1/2})}/{(q-1)} $
which leads to
\begin{align*}
	c_n(q)=\frac{q^{n+\frac{3}{2}}}{1-q}h_{n+1}.
\end{align*}
As a consequence, we have the following evaluation of the moment Pfaffian
\begin{align*}
	\Pf(m_{i,j})_{i,j=0}^{2N-1}=q^{-\frac{1}{6}N(N+1)(8N-5)}\prod_{k=0}^{N-1}(q;q)_{2k+1}.
\end{align*}
On the other hand, we can compute the skew moments \eqref{sm} as a scalar product of canonical moments $ \mu_{n} $
\begin{align*}
		m_{i,j}=&([j]_q-[i]_q)\int_{0}^{\infty}x^{i+j-1}\omega(x;q)\dq=([j]_q-[i]_q)\int_{0}^{\infty}x^{i+j+1}q^{\frac{3}{2}}\rho(x;q)\dq\\ =&([j]_q-[i]_q)q^{{3}/{2}}\mu_{i+j+1}=([j]_q-[i]_q)q^{-[(i+j+2)^2-3]/2}.
\end{align*}
By combining above results, we have the evaluation
\begin{align*}
	\Pf(([j]_q-[i]_q)q^{-(i+j+2)^2/2})_{i,j=0}^{2N-1}=q^{-\frac{1}{6}N(2N+1)(8N-1)}\prod_{k=0}^{N-1}(q;q)_{2k+1}.
\end{align*}

\subsection{Little $q$-Jacobi case}
Little $q$-Jacobi polynomials are important in many mathematical fields such as polynomials theory \cite{koekoek10} and quantum group \cite{masuda91}.
These polynomials have the following series form (c.f. \cite[eq. (2.21)]{masuda91})
\begin{align*}
p_n^{(\alpha,\beta)}(z;q)=\sum_{r\geq0}\frac{(q^{-n};q)_r(q^{\alpha+\beta+n+1};q)_r}{(q;q)_r(q^{\alpha+1};q)_r}(qz)^r
\end{align*}
and they obey the following orthogonal relation \cite[Prop. 3.9]{masuda91}
\begin{align*}
\int_0^1 p_m^{(\alpha,\beta)}(z;q)p_n^{(\alpha,\beta)}(z;q)z^\alpha(qz;q)_\beta d_qz=\delta_{n,m}q^{(\alpha+1)n}\frac{(1-q)(q;q)_\alpha^2(q;q)_{\beta+n}(q;q)_n}{(1-q^{\alpha+\be+2n+1})(q;q)_{\alpha+n}(q;q)_{\alpha+\beta+n}}.
\end{align*}
In this case, the canonical moments are given by \cite{njionou}
\begin{align*}
	\mu_{n}^{(\alpha,\beta)}=\int_0^1 z^{\al+n}(qz;q)_\be d_qz=\frac{(1-q)(q^{\alpha+\beta+1};q)_\infty(q;q)_\infty(q^{\alpha+1};q)_n}{(1-q^{\alpha+\beta+1})(q^{\alpha+1};q)_\infty(q^{\beta+1};q)_\infty(q^{\alpha+\beta+2};q)_n}.
\end{align*}
If we define monic little $q$-Jacobi polynomials $\{\tilde{p}_n^{(\al,\be)}(z;q)\}_{n\geq0}$
\begin{align*}
\tilde{p}_n^{(\al,\be)}(z;q)=(-1)^n q^{\frac{n(n-1)}{2}}\frac{(q^{\al+1};q)_n}{(q^{\al+\be+n+1};q)_n}p_n^{(\al,\be)}(z;q),
\end{align*}
then they satisfy the orthogonal relation
\begin{align}\label{lqj}
\begin{aligned}
\int_0^1 &\tilde{p}_n^{(\al,\be)}(z;q)\tilde{p}_m^{(\al,\be)}(z;q)z^\al (qz;q)_\be d_qz\\
&=q^{n(n+\al)}[\al+\be+2n+1]_q^{-1}\frac{(q;q)_n(q;q)_{n+\al}(q;q)_{n+\be}(q;q)_{n+\al+\be}}{(q;q)_{2n+\al+\be}^2}\delta_{n,m}:=h_n\delta_{n,m}.
\end{aligned}
\end{align}
Regarding the weight of little $q$-Jacobi polynomials
\begin{align*}
\rho(z;q)=z^\al (qz;q)_\be,
\end{align*} 
by solving \eqref{pearson}, it admits the Pearson pair
\begin{align*}
	(f,g)=\left(-x^{2}+x,-q^{\frac{1}{2}}\left([\alpha+\beta+2]_{q} x-[\alpha+1]_{q}\right)\right).
\end{align*}
Therefore, the coefficient $ c_n $ has the expression
\begin{align*}
	c_n=q^{-n}[2n+2+\alpha+\beta]_qh_{n+1},
\end{align*}
where $ h_n $ is the nomalization constant in the orthogonal relation of monic polynomials \eqref{lqj}.
Then we have the following expression for the moment Pfaffian
\begin{align*}
	\Pf(m_{ij}^{(\alpha,\beta)})_{i,j=0}^{2N-1}&=q^{\frac{1}{3}N(4N^2-3N+2)+\alpha N^2}\\
	&\times\prod_{k=0}^{N-1}\frac{(1-q^{\alpha+\beta+4k+2})(q,q^{\alpha+1},q^{\beta+1},q^{\alpha+\beta+1};q)_{2k+1}(q;q)_\infty(q^{\alpha+\beta+1};q)_\infty}{(q^{\alpha+1},q^{\beta+1};q)_\infty(q^{\alpha+\beta+1};q)_{4k+2}(q^{\alpha+\beta+1};q)_{4k+3}}.
\end{align*}
Since the weight of the corresponding skew orthogonal little $q$-Jacobi polynomials is 
\begin{align*}
\omega(x;q)=f(qx;q)\rho^{(\al,\be)}(qx;q)=q^{\alpha+1}\rho^{(\alpha+1,\beta+1)}(x;q)=q^{\alpha+1}x^{\alpha+1}(qx;q)_{\beta+1},
\end{align*}
 the skew moments $\{ m_{i,j}^{(\alpha,\beta)} \}_{i,j\in\mathbb{N}}$ are related to the canonical moments by
\begin{align*}
	m_{i,j}^{(\alpha,\beta)}=&([j]_q-[i]_q)\int_{0}^{1}x^{i+j-1}\omega(x;q)\dq=([j]_q-[i]_q)q^{\alpha+1}\mu_{i+j-1}^{(\alpha+1,\beta+1)}\\
	=&([j]_q-[i]_q)q^{\alpha+1}\frac{(1-q)(q^{\alpha+\beta+1};q)_\infty(q;q)_\infty(q^{\alpha+1};q)_n}{(1-q^{\alpha+\beta+1})(q^{\alpha+1};q)_\infty(q^{\beta+1};q)_\infty(q^{\alpha+\beta+2};q)_n}.
\end{align*}
After eliminating the constant $ \left[q^{\alpha+1}{(q^{\alpha+\beta+1};q)_\infty(q;q)_\infty}/{\left((1-q^{\alpha+\beta+1})(q^{\alpha+1},q^{\beta+1};q)_\infty\right)}\right]^N $, we have
\begin{align}\label{lj}
&\Pf\left(([j]_q-[i]_q)\frac{(1-q)(q^{\al+2};q)_{i+j-1}}{(q^{\al+\be+4};q)_{i+j-1}}\right)_{i,j=0}^{2N-1}\\
&\qquad=q^{\frac{1}{3}N(N-1)(4N+1)+\al N(N-1)}\prod_{k=0}^{N-1}\frac{(q;q)_{2k+1}(q^{
\al+2},q^{\be+2};q)_{2k}}{(q^{\al+\be+4};q)_{2k-2}(q^{\al+\be+2k+2};q)_{2k}(q^{\al+\be+2k+2};q)_{2k+2}},\nonumber
\end{align}
which coincides with the result in \cite[Thm. 5.1]{ishikawa20}.

\subsection{Big $q$-Jacobi case}
Big $ q $-Jacobi polynomials were introduced by Andrews and Askey as an infinite-dimentional version of the $ q $-Hahn polynomials \cite{askey}. In addition, big $q$-Jacobi polynomials were also contained in the Bannai-Ito scheme of
dual systems of orthogonal polynomials as an infinite dimension analogue of the q-Racah polynomials \cite{Bannai}. These polynomials take the hypergeometric function form
\begin{align*}
P_{n}(x ; a, b, c ; q)={ }_{3} \phi_{2}\left(\left.\begin{array}{c}
q^{-n}, a b q^{n+1}, x \\
a q, c q
\end{array} \right| q ; q\right),
\end{align*}
orthogonal with respect to the weight
\begin{align*}
\rho^{(a,b,c)}(x;q)=\frac{(a^{-1}x,c^{-1}x;q)_\infty}{(x,bc^{-1}x;q)_\infty}
\end{align*}
and have the orthogonality
\begin{align*}
&\int_{c q}^{a q} P_{m}(x ; a, b, c ; q) P_{n}(x ; a, b, c ; q) \rho^{(a,b,c)}(x;q)d_{q} x \\
&\quad=a q(1-q) \frac{\left(q, a^{-1} c, a c^{-1} q, a b q^{2} ; q\right)_{\infty}}{\left(a q, b q, c q, a b c^{-1} q ; q\right)_{\infty}}
 \frac{(1-a b q)}{\left(1-a b q^{2 n+1}\right)} \frac{\left(q, b q, a b c^{-1} q ; q\right)_{n}}{(a b q, a q, c q ; q)_{n}}\left(-a c q^{2}\right)^{n} q^{\binom{n}{2}} \delta_{m n}.
\end{align*}
One can easily check that the normalisation constant for the monic polynomials is
\begin{align*}
h_n=\frac{aq(1-q)(-acq^2)^nq^{\binom{n}{2}}(q,a^{-1}c,ac^{-1}q,abq;q)_\infty(q,aq,bq,cq,abc^{-1}q;q)_n}{(1-abq^{2n+1})(aq,bq,cq,abc^{-1}q;q)_\infty(abq;q)_n(abq^{n+1};q)_n^2}
\end{align*}
and the canonical moments are given by \cite{njionou}
\begin{align*}
	\mu_{n}^{(a,b,c)}&=\int_{cq}^{aq}x^n \rho^{(a,b,c)}(x;q)d_qx\\&=a q \frac{\left(a b q^{2}, a^{-1} c, a c^{-1} q ; q\right)_{\infty}}{\left(a q, b q, c q, a b c^{-1} q ; q\right)_{\infty}} \sum_{m=0}^{n}(-1)^{m}\left[\begin{array}{c}
	n \\
	m
	\end{array}\right]_{q} q^{-n m+\binom{m+1}{2}} \frac{(a q, c q ; q)_{m}}{\left(a b q^{2} ; q\right)_{m}}.
\end{align*}
Solving equation (\ref{pearson}) gives the Pearson pair
\begin{align*}
(f,g)=\left((1-\frac{x}{aq})(1-\frac{x}{cq}),\frac{q^{\frac{1}{2}}}{1-q}\left((\frac{1}{acq^2}-\frac{b}{c})x+\frac{b}{c}+1-\frac{1}{aq}-\frac{1}{cq}\right)\right).
\end{align*}
Then, by comparing the leading coefficients on both sides of (\ref{eq:cn}) we have
\begin{align*}
c_n=\frac{abq^{2n+2}-1}{ac(1-q)q^{n+2}}h_{n+1},
\end{align*}
which leads to the following evaluation of the moment Pfaffian
\begin{align*}
\Pf(m_{ij}^{(a,b,c)})_{i,j=0}^{2N-1}=&(-1)^Nc^{N(N-1)}a^{N^2}q^{\frac{1}{6}N(N+1)(4N-1)}\\
\times&\prod_{k=0}^{N-1}\frac{(q,a^{-1}c,ac^{-1}q,abq;q)_\infty(q,aq,bq,cq,abc^{-1}q;q)_{2k+1}}{(aq,bq,cq,abc^{-1}q;q)_\infty(abq;q)_{2(k+N)+1}}.
\end{align*}
where $m_{i,j}^{(a,b,c)}$ is given by
\begin{align*}
m_{i,j}^{(a,b,c)}=([j]_q-[i]_q)\int_{cq}^{aq} x^{i+j-1}\omega(x;q)^{(a,b,c)}d_qx,
\end{align*}
and $\omega(x;q)$ is expressed by
\begin{align*}
 \omega(x;q)^{(a,b,c)}=f(qx;q)\rho^{(a,b,c)}(qx;q)=(1-x)(1-bc^{-1}x)\rho(x;q)^{(a,b,c)}=\rho^{(qa,qb,qc)}(qx;q).
 \end{align*}
Note that
\begin{align*}
\mu_{n}^{(a,b,c)}=\int_{cq}^{aq}x^n\rho^{(a,b,c)}(x;q)\dq=q^{n+1}\int_{c}^{a}x^n\rho^{(a,b,c)}(qx;q)\dq,
\end{align*}
we have the expression
\begin{align*}
m^{(a,b,c)}_{i,j}=&([j]_q-[i]_q)q^{-(i+j)}\mu_{i+j-1}^{(qa,qb,qc)}\\
=&([j]_q-[i]_q)aq^{-(i+j-1)}\frac{\left(a b q^{2}, a^{-1} c, a c^{-1} q ; q\right)_{\infty}}{\left(a q, b q, c q, a b c^{-1} q ; q\right)_{\infty}}\left(\sum_{m=0}^{n}(-1)^{m}\left[\begin{array}{c}
n \\
m
\end{array}\right]_{q} q^{-n m+\binom{m+1}{2}} \frac{(a q, c q ; q)_{m}}{\left(a b q^{2} ; q\right)_{m}}\right).
\end{align*}
Dividing both sides by $ \left[aq{\left(a b q^{2}, a^{-1} c, a c^{-1} q ; q\right)_{\infty}}/{\left(a q, b q, c q, a b c^{-1} q ; q\right)_{\infty}}\right]^N $ and changing variables $ (a,b,c)\rightarrow(q^{-1}a,q^{-1}b,q^{-1}c) $, we get
\begin{align*}
\Pf&\left( ([j]_q-[i]_q)\sum_{m=0}^{i+j-1}(-1)^{m}\left[\begin{array}{c}
i+j-1 \\
m
\end{array}\right]_{q} q^{-(i+j-1)(m+1)+\binom{m+1}{2}} \frac{(a q, c q ; q)_{m}}{\left(a b q^{2} ; q\right)_{m}} \right)_{i,j=0}^{2N-1}\\
&=(-1)^N(ac)^{N(N-1)}q^{\frac{1}{6}N(4N^2-9N+11)}\prod^{N-1}_{k=0}\frac{(q;q)_{\infty}(q;q)_{2k+1}(a,b,c,abc^{-1}q;q)_{2k}}{(abq^2;q)_{2(k+N)-2}}.
\end{align*}
\begin{remark}
With $ a=b=1 $, big q-Jacobi polynomials reduce to the big q-Legendre polynomials \cite{li20}. In particular, we have the following evaluation of q-Catalan-Hankel Pfaffian
\begin{align*}
	\Pf&\left( ([j]_q-[i]_q)\sum_{m=0}^{i+j-1}(-1)^{m}\left[\begin{array}{c}
	i+j-1 \\
	m
	\end{array}\right]_{q} q^{-(i+j-1)(m+1)+\binom{m+1}{2}} \frac{( q, c q ; q)_{m}}{\left( q^{2} ; q\right)_{m}} \right)_{i,j=0}^{2N-1}\\
	&=(-1)^N(c)^{N(N-1)}q^{\frac{1}{6}N(4N^2-9N+11)}\prod^{N-1}_{k=0}\frac{(q;q)_{\infty}(q;q)_{2k+1}(q,q,cq,c^{-1}q;q)_{2k}}{(q^2;q)_{2(k+N)-2}}.
\end{align*}
	
\end{remark}

\section{Further remarks}
In this paper, we have developed a method based on the relation between classical ($ q $-)orthogonal and ($ q $-)skew orthogonal polynomials to evaluate certain $q$-Catalan-Hankel Pfaffians whose entries are composed of the moments of classical orthogonal polynomials. Some examples are given to illustrate the approach including the continuous ones (e.g. Hermite, Laguerre, Jacobi and Cauchy) and discrete $q$-cases (e.g. Al-Salam $\&$ Carlitz I, little $q$-Jacobi, Stieltjes-Wigert and big $q$-Jacobi polynomials). Amoung those examples, the Al-Salam $\&$ Carlitz I and Little $q$-Jacobi case are compared with the results obtained by M. Ishikawa and J. Zeng in \cite{ishikawa20} and we present alternative proofs of \cite[eq. (6.7) \& Thm 5.2]{ishikawa20}. Besides, the examples in \cite[Conjecture 7.1]{ishikawa20} seem to be related to some discrete measure on the linear lattice. However, as mentioned in \cite{forrester20}, the skew moments defined in linear lattices are of the form $\mu_{i,j}=(j-i)\mu_{i+j-1}+({j\choose2}-{i\choose 2})\mu_{i+j-2}+\cdots$, it is unclear to us whether it could be written in the Catalan-Hankel Pfaffian form.

As mentioned before, $ q $-analogues of Selberg integral can be used to evaluate $q$-Catalan-Hankel Pfaffians. For example, the formula \eqref{lj} can be evaluated by the following Askey-Habsieger-Kadell formula \cite{ishikawa20}
\begin{align*}
	&\int_{[0,1]^{n}} \prod_{i<j} t_{i}^{2 k}\left(q^{1-k} t_{j} / t_{i} ; q\right)_{2 k} \prod_{i=1}^{n} t_{i}^{x-1} \frac{\left(t_{i} q ; q\right)_{\infty}}{\left(t_{i} q^{y} ; q\right)_{\infty}} d_{q} \boldsymbol{t}\\
	&\quad=q^{k x\binom{n}{2}+2 k^{2}\binom{n}{3}} \prod_{j=1}^{n} \frac{\Gamma_{q}(x+(j-1) k) \Gamma_{q}(y+(j-1) k) \Gamma_{q}(j k+1)}{\Gamma_{q}(x+y+(n+j-2) k) \Gamma_{q}(k+1)}.
\end{align*}
Conversely, in this paper, we give an explanation to this integral identity in terms of the skew little $ q $-Jacobi polynomials.
It is still an open question to construct classical skew orthogonal polynomials to make explanations of some other discrete Selberg integrals such as \cite[Prop. 5.1]{brent16}
\begin{align*}
&\sum_{k_1,\cdots,k_r=-n}^n \prod_{1\leq i<j\leq r}|[k_i-k_j]_q|\prod_{i=1}^rq^{(k_i+n-r+i)^2/2}\left[
\begin{array}{c}
2n\\n+k_i\end{array}
\right]\\
&\quad=\frac{r!}{|r|_{q^{1/2}}!}\prod_{i=1}^r (-q^{1/2};q^{1/2})_i(-q^{i/2+1};q)_{2n-r}\prod_{i=1}^r\frac{\Gamma_q(1+i/2)}{\Gamma_q(3/2)}\frac{\Gamma_q(2n+1)\Gamma_q(2n-i+5/2)}{\Gamma_q(2n-i+2)\Gamma_q(2n-i/2+2)}.
\end{align*}
Moreover, combinatoric explanations of $q$-Catalan-Hankel Pfaffians are left to be further investigated.

\small
\bibliographystyle{abbrv}

\end{document}